\magnification=1202

\input amssym.def
\font\twlveusm=eusm10 at9.98pt
\font\seveneusm=eusm7 at6.38pt 
\font\fiveeusm=eusm5
\newfam\eusmfam
\textfont\eusmfam=\twlveusm
\scriptfont\eusmfam=\seveneusm
\scriptscriptfont\eusmfam=\fiveeusm
\def\eusm#1{{\fam\eusmfam\relax#1}}
 at9.98pt
\font\bb=msbm10 at9.98pt
 at9.98pt 
 at5pt
\font\cyr=wncyi10 at9.98pt
\font\eightrm=cmr8
\font\eightsc=cmcsc8
\font\labf=cmbx10 at13.1pt
 at15.74pt
\font\larm=cmr10 at13.1pt
 at15.74pt
\font\sc=cmcsc10 at9.98pt
\font\tenpbf=cmbx10 at8.32pt
\font\tenpit=cmti10 at8.32pt
\font\tenprm=cmr10 at8.32pt


\def\Ad{{\rm Ad}}

\def\an{\raise0.5pt\hbox{$\kern2pt\scriptstyle\in\kern2pt$}}
\def\Ann{\hbox{\rm Ann\kern1pt}}
\def\Anns{\hbox{$\scriptstyle\rm Ann\kern0.5pt$}}

\def\arkef{\advance\chapternumber by 1\sc\roman{\the\chapternumber}}
\def\Aut{\hbox{\rm Aut\kern1pt}}
\def\bell{\hskip0pt\lower1.6pt\hbox{\bel\char'012}\kern5pt}

\def\callige{\hbox{\calligl e\kern2pt}}
\def\cheridexi{\hskip0pt\lower2pt\hbox{\cheridexia}\kern5pt}
\def\cheridexia{{\bbding\char'21}}
\def\coad{{\rm coad}}
\def\Coad{{\rm Coad}}
\def\coker{{\rm coker\kern1pt}}
\def\Colon{\colon\kern2pt}
\def\comp{\hbox{\lower5.8pt\hbox{\larm\char'027}}}

\def\corang{{\rm corang\kern1pt}}
\def\cos{\hbox{\rm cos\kern1pt}}
\def\cosh{\hbox{\rm cosh\kern1pt}}
\def\dbaraux{\hbox{\= {\kern-2pt\= {}}}}
\def\dbar#1{\raise3pt\hbox{\dbaraux}\kern-7.8pt #1}
\def\Der{\lower0.5pt\hbox{\ygoth Der}}

\def\dim{{\rm dim\kern1pt}}
\def\double{\hbox{\kern1.5pt\bb\char'156\kern-7.6pt\char'157\kern1.5pt}}

\def\enwsh#1{{\lower2.1pt\hbox{$\buildrel{\textstyle\cup}\over
{\lower.8pt\hbox{${}_{\scriptscriptstyle#1}$}}$}}}
\def\exp{{\rm exp\kern1pt}}

\def\exten{\hbox{\callig \kern-2.5pt Ext\lower2.5pt\hbox{\kern2.5pt}}}
\def\Ham{\hbox{\rm Ham\kern1pt}}
\def\im{{\rm im\kern1pt}}
\def\k{\raise0.25pt\hbox{$\ygot k$}}

\def\ker{{\rm ker\kern1pt}}

\def\Lie{\hbox{
\callig Lie\kern2pt}}
\def\mavrodexi{\hskip0pt\lower2pt\hbox{\mavrodexia}\kern5pt}
\def\mavrodexia{{\bbding\char'15}}
\def\meriki{\hbox{\cyr\char'144\kern0.3pt}}

\def\na{\raise0.5pt\hbox{$\kern2pt\scriptstyle\ni\kern2pt$}}
\def\noan{\hbox{$\an\raise0.6pt\hbox{$\kern-6.5pt\scriptstyle
          \slash\kern3pt$}$}}
\def\odot{\;{\mathchar"220C}\;}

\def\oplus{\;{\mathchar"2208}\;}

\def\pounds{\rlap{\lower3.5pt\hbox{\kern2.9pt\hbox{\char'26}}}
           {\script L}}
\def\pr{\hbox{\kern3pt{\calligs p}\callig r\kern2pt}}
\def\qed{\hbox{\kern0.3cm\vrule height5pt width5pt depth-0.2pt}}
\def\QED{\hbox{\kern0.3cm\vrule height6pt width6pt depth-0.2pt}}
\def\R{\hbox{\bf\char'122}}
\def\rang{{\rm rang\kern1pt}}
\def\rank{{\rm rank\kern1pt}}
\def\S{\hbox{\bf\char'123}}
\def\san{\raise0.5pt\hbox{$\kern0.7pt\scriptscriptstyle
         \in\kern0.7pt$}}

\def\scomp{\hskip-0.05truecm\hbox{\lower5pt\hbox{$\mathchar"2017$}}
           \hskip-0.05truecm}

\def\sem{\hbox{{\script S}\kern-2.5pt\callig em\kern2pt}}

\def\sin{\hbox{\rm sin\kern1pt}}
\def\sinh{\hbox{\rm sinh\kern1pt}}

\def\styl{\hbox{\bbding\char'26}}
\def\stylo{\hskip0.3truecm\hbox{\lower1.5pt\hbox{\styl}}}
\def\times{\;{\mathchar"2202}\;}
\def\timess#1{\hbox{$\times_{\hskip-0.1truecm#1\hskip0.14truecm}$}}
\def\tonos{\hbox{\kern-1.3pt\lower0.7pt\hbox{$\mathchar"6013$}}}
\def\tonoskef{\hbox{$\kern-1.3pt\mathchar"6013$}}

\def\wbaraux{\hbox{\= {\kern-1.4pt\= {\kern-1.4pt\= {\kern-1.4pt\=
 {\kern-1.4pt\= {\kern-1.4pt\= {\kern-1.4pt\= {\kern-1.4pt\= {}}}}}}}}}}
\def\wbar#1{\hbox{\raise3pt\hbox{\wbaraux}\kern-30.5pt #1}}

\def\wwbaraux{\hbox{\= {\kern-1.4pt\= {\kern-1.4pt\= {\kern-1.4pt\=
{\kern-1.4pt\= {\kern-1.4pt\= {\kern-1.4pt\= {\kern-1.4pt\=
{\kern-1.4pt\= {}}}}}}}}}}}
\def\wwbar#1{\hbox{\raise3pt\hbox{\wwbaraux}\kern-34pt #1}}

\catcode`\@=11
\def\eightpoint{\eightrm}
\def\footnote#1{\edef\@sf{\spacefactor\the\spacefactor}#1\@sf
     \insert\footins\bgroup\eightpoint
     \interlinepenalty100 \let\par=\endgraf
      \leftskip=0pt \rightskip=0pt
      \splittopskip=10pt plus 1pt minus 1pt \floatingpenalty=20000
      \smallskip\item{#1}\bgroup\strut\aftergroup\@foot\let\neft}
\skip\footins=12pt plus 2pt minus 4pt
\dimen\footins=30pc

\def\line{\hbox to\hsize}

\def\title#1{\line{\hss}\line{\hss#1\hss}%
\line{\hss}\hskip-0.75truecm}

\def\author#1{{\tenprm #1:}}
\def\ekdoths#1{{\tenprm #1}}

\def\periodiko#1{{\tenpit #1\tenprm ,}}
\def\selides#1{{\tenprm #1}}
\def\titlosa#1{{\tenprm #1,}}
\def\titlosb#1{{\tenpit #1\tenprm ,}}
\def\volume#1{{\tenprm Vol. \tenpbf #1\tenprm :}}

%
%
%
\def\teleia{\hbox{.}}
\newif\ifPhysRev
\def\Textindent#1{\noindent\llap{#1\enspace}\ignorespaces}
\def\nonfrenchspacing{\sfcode`\.=3001 \sfcode`\!=3000 \sfcode`\?=3000
        \sfcode`\:=2000 \sfcode`\;=1500 \sfcode`\,=1251 }
\nonfrenchspacing
\newdimen\d@twidth
 {\setbox0=\hbox{s.} \global\d@twidth=\wd0 \setbox0=\hbox{s}
        \global\advance\d@twidth by -\wd0 }
\def\removehglue{\loop \unskip \ifdim\lastskip >\z@ \repeat }
\def\roll@ver#1{\removehglue \nobreak \count255 =\spacefactor \dimen@=\z@
        \ifnum\count255 =3001 \dimen@=\d@twidth \fi
        \ifnum\count255 =1251 \dimen@=\d@twidth \fi
    \iftwelv@ \kern-\dimen@ \else \kern-0.83\dimen@ \fi
   #1\spacefactor=\count255 }
\def\step@ver#1{\relax \ifmmode #1\else \ifhmode
        \roll@ver{${}#1$}\else {\setbox0=\hbox{${}#1$}}\fi\fi }
\def\attach#1{\step@ver{\strut^{\mkern 2mu #1} }}

\normalbaselineskip = 20pt plus 0.2pt minus 0.1pt
\normallineskip = 1.5pt plus 0.1pt minus 0.1pt
\normallineskiplimit = 1.5pt
\newskip\normaldisplayskip
\normaldisplayskip = 20pt plus 5pt minus 10pt
\newskip\normaldispshortskip
\normaldispshortskip = 6pt plus 5pt
\newskip\normalparskip
\normalparskip = 6pt plus 2pt minus 1pt
\newskip\skipregister
\skipregister = 5pt plus 2pt minus 1.5pt
\newif\ifsingl@    \newif\ifdoubl@
\newif\iftwelv@    \twelv@true
\def\singlespace{\singl@true\doubl@false\spaces@t}
\def\doublespace{\singl@false\doubl@true\spaces@t}
\def\normalspace{\singl@false\doubl@false\spaces@t}
\def\Tenpoint{\tenpoint\twelv@false\spaces@t}
\def\Twelvepoint{\twelvepoint\twelv@true\spaces@t}
\def\spaces@t{\relax
      \iftwelv@ \ifsingl@\subspaces@t3:4;\else\subspaces@t1:1;\fi
       \else \ifsingl@\subspaces@t3:5;\else\subspaces@t4:5;\fi \fi
      \ifdoubl@ \multiply\baselineskip by 5
         \divide\baselineskip by 4 \fi }
\def\subspaces@t#1:#2;{
      \baselineskip = \normalbaselineskip
      \multiply\baselineskip by #1 \divide\baselineskip by #2
      \lineskip = \normallineskip
      \multiply\lineskip by #1 \divide\lineskip by #2
      \lineskiplimit = \normallineskiplimit
      \multiply\lineskiplimit by #1 \divide\lineskiplimit by #2
      \parskip = \normalparskip
      \multiply\parskip by #1 \divide\parskip by #2
      \abovedisplayskip = \normaldisplayskip
      \multiply\abovedisplayskip by #1 \divide\abovedisplayskip by #2
      \belowdisplayskip = \abovedisplayskip
      \abovedisplayshortskip = \normaldispshortskip
      \multiply\abovedisplayshortskip by #1
        \divide\abovedisplayshortskip by #2
      \belowdisplayshortskip = \abovedisplayshortskip
      \advance\belowdisplayshortskip by \belowdisplayskip
      \divide\belowdisplayshortskip by 2
      \smallskipamount = \skipregister
      \multiply\smallskipamount by #1 \divide\smallskipamount by #2
      \medskipamount = \smallskipamount \multiply\medskipamount by 2
      \bigskipamount = \smallskipamount \multiply\bigskipamount by 4 }
\def\normalbaselines{ \baselineskip=\normalbaselineskip
   \lineskip=\normallineskip \lineskiplimit=\normallineskip
   \iftwelv@\else \multiply\baselineskip by 4 \divide\baselineskip by 5
     \multiply\lineskiplimit by 4 \divide\lineskiplimit by 5
     \multiply\lineskip by 4 \divide\lineskip by 5 \fi }


\def\abstract#1{\parshape=1 0.7cm \dimen10
                {\tenpbf Abstract. \tenprm #1}}

\newcount\appendixnumber     \appendixnumber=0
\newcount\chapternumber      \chapternumber=0
\newcount\equanumber         \equanumber=0
\newcount\mathnumber         \mathnumber=0
\newcount\appequanumber      \appequanumber=0
\newcount\appmathnumber      \appmathnumber=0

\let\variableone=\relax
\let\variabletwo=\relax
\let\chapterlabel=\relax
\let\sectionlabel=\relax
\let\mathlabel=\relax
\newtoks\chapterstyle        \chapterstyle={\Number}
\newtoks\sectionstyle        \sectionstyle={\chapterlabel\Number}
\newskip\chapterskip         \chapterskip=\bigskipamount
\newskip\sectionskip         \sectionskip=\medskipamount
\newskip\headskip            \headskip=8pt plus 3pt minus 3pt
\newdimen\chapterminspace    \chapterminspace=15pc
\newdimen\sectionminspace    \sectionminspace=10pc
\newdimen\sectionspace       \sectionspace=20pc
\newdimen\referenceminspace  \referenceminspace=25pc

\def\chapterreset{\global\advance\chapternumber by 1
   \ifnum\equanumber<0 \else\global\equanumber=0\fi
   \mathnumber=0
   \makechapterlabel}
\def\makechapterlabel{\let\sectionlabel=\relax\let\mathlabel=\relax
 \xdef\chapterlabel{\the\chapterstyle{\the\chapternumber\teleia\kern3pt}}}

\def\rightheadline{\sc\hfil\variableone\eightsc\hfil\folio}
\def\leftheadline{\eightsc\folio\hfil{\sc\variabletwo}\hfil}
\def\heads{\footline={\hfil}\headline={\ifodd\pageno
               \rightheadline\else\leftheadline\fi}}

\def\headseis{\partreset\headline={\ifodd\pageno{
                         \hfil\sc partie {\eightsc\the\partnumber}
                         -introduction\hfil\eightsc\folio}\else
                        {\eightsc\folio\hfil\sc partie
                         {\eightsc\the\partnumber}-introduction\hfil}\fi}
                        \footline={\hfil}}

\def\alphabetic#1{\count255='140 \advance\count255 by #1\char\count255}
\def\Alphabetic#1{\count255='100 \advance\count255 by #1\char\count255}
\def\Roman#1{\uppercase\expandafter{\romannumeral #1}}
\def\roman#1{\romannumeral #1}
\def\Number#1{\number #1}
\def\BLANC#1{}

\def\titlestyle#1{\par\begingroup \interlinepenalty=9999
     \leftskip=0.02\hsize plus 0.23\hsize minus 0.02\hsize
     \rightskip=\leftskip \parfillskip=0pt
     \hyphenpenalty=9000 \exhyphenpenalty=9000
     \tolerance=9999 \pretolerance=9000
     \spaceskip=0.333em \xspaceskip=0.5em
     \iftwelv@\bf\else\bf\fi
   \noindent #1\par\endgroup }

\def\spacecheck#1{\dimen@=\pagegoal\advance\dimen@ by -\pagetotal
   \ifdim\dimen@<#1 \ifdim\dimen@>0pt \vfil\break \fi\fi}
\def\TableOfContentEntry#1#2#3{\relax}

\def\chapter#1{\par\vskip0.7cm
   \chapterreset \titlestyle{\chapterlabel\ #1}
   \nobreak\vskip\headskip
   \wlog{\string\chapter\space \chapterlabel} }

\def\appendixreset{\global\advance\appendixnumber by 1
                   \appmathnumber=0\appequanumber=0}
\def\appendix#1{\par \penalty-300\vskip\chapterskip
   \spacecheck\chapterminspace
   \appendixreset \title{\bf Appendix \Alphabetic{\the\appendixnumber}}
   \nobreak\vskip-\chapterskip\penalty 30000
   \vskip-\chapterskip
   \par{\titlestyle{#1}}
   \vskip\chapterskip
   \wlog{\string\appendix\space \chapterlabel} }

%
%
\def\eqname#1{\relax \ifnum\equanumber<0
     \xdef#1{{\noexpand\rm(\number-\equanumber)}}%
       \global\advance\equanumber by -1
    \else \global\advance\equanumber by 1
      \xdef#1{{\noexpand(\rm{\the\chapternumber}\teleia
                            \rm{\number\equanumber})}} \fi #1}

\def\eqn{\eqno\eqname}

\def\math#1#2{\vskip0.1cm
   \global\advance\mathnumber by 1
   \xdef\mathlabel{\the\chapternumber\teleia\the\mathnumber}
   \wlog{\string\math\space \mathlabel}
   {\bf\enspace\mathlabel\hskip0.2cm #1}
   \xdef#2{{\mathlabel}}}

\def\appeqname#1{\relax \ifnum\appequanumber<0
     \xdef#1{{\noexpand\rm(\number-\appequanumber)}}%
       \global\advance\appequanumber by -1
    \else \global\advance\appequanumber by 1
      \xdef#1{{\noexpand(\hbox{\Alphabetic{\the\appendixnumber}}\teleia
                            {\number\appequanumber})}} \fi #1}

\def\mathapp#1#2{\vskip0.1cm
   \global\advance\appmathnumber by 1
   \xdef\appmathlabel{{\Alphabetic{\the\appendixnumber}}\teleia
   \the\appmathnumber}
   \wlog{\string\mathapp\space \appmathlabel}
   {\bf\enspace\appmathlabel\hskip0.2cm #1}
   \xdef#2{{\appmathlabel}}}


%
%
%
\newtoks\referencestyle      \referencestyle={\tenpbf\Number}
\newcount\referencecount     \referencecount=0
\newcount\lastrefsbegincount \lastrefsbegincount=0
\newif\ifreferenceopen       \newwrite\referencewrite
\newif\ifrw@trailer
\newdimen\refindent     \refindent=13pt
\def\NPrefmark#1{\attach{\scriptscriptstyle [ #1 ] }}
\let\PRrefmark=\attach
\def\refmark#1{\relax\ifPhysRev\PRrefmark{#1}\else\NPrefmark{#1}\fi}
\def\refend@{\refmark{\number\referencecount}}
\def\refend{\refend@{}\space }
\def\refsend{\refmark{\count255=\referencecount
   \advance\count255 by-\lastrefsbegincount
   \ifcase\count255 \number\referencecount
   \or \number\lastrefsbegincount,\number\referencecount
   \else \number\lastrefsbegincount-\number\referencecount \fi}\space }
\def\refitem#1{\par\hangafter=0 \hangindent=\refindent	\Textindent{#1}}
\def\Ref{\rw@trailertrue\REF}
\def\REF#1{\r@fstart{#1}%
   \rw@begin{\tenprm [\tenpbf\Number{\the\referencecount}\tenprm ]}\rw@end}
\def\r@fstart#1{\chardef\rw@write=\referencewrite \let\rw@ending=\refend@
   \ifreferenceopen \else \global\referenceopentrue
   \immediate\openout\referencewrite=referenc.txa
   \toks0={\catcode`\^^M=10}\immediate\write\rw@write{\the\toks0} \fi
   \global\advance\referencecount by 1 
   \xdef#1{[{\the\referencestyle{\the\referencecount}}]}}
 {\catcode`\^^M=\active %
 \gdef\rw@begin#1{\immediate\write\rw@write{\noexpand\refitem{#1}}%
   \begingroup \catcode`\^^M=\active \let^^M=\relax}%
 \gdef\rw@end#1{\rw@@end #1^^M\rw@terminate \endgroup%
   \ifrw@trailer\rw@ending\global\rw@trailerfalse\fi }%
 \gdef\rw@@end#1^^M{\toks0={#1}\immediate\write\rw@write{\the\toks0}%
   \futurelet\n@xt\rw@test}%
 \gdef\rw@test{\ifx\n@xt\rw@terminate \let\n@xt=\relax%
       \else \let\n@xt=\rw@@end \fi \n@xt}%
}
\let\rw@ending=\relax
\let\rw@terminate=\relax

\def\Closeout#1{\toks0={\catcode`\^^M=5}\immediate\write#1{\the\toks0}%
   \immediate\closeout#1}

\topskip1truecm
\voffset=2.5truecm
\hsize 15truecm
\vsize 20truecm
\hoffset=0.5truecm
\def\undertext#1{$\underline{\hbox{#1}}$}
\topglue 3truecm
\nopagenumbers
\def\variableone{p. baguis}
\def\variabletwo{homogeneous symplectic manifolds of poisson-lie 
groups}
\def\absize{11cm}

\def\abstract#1{\baselineskip=12pt plus .2pt
                \parshape=1 0.7cm \absize
                 {\tenpbf Abstract. \tenprm #1}}



\Ref\abmar{
\author{Abraham, R., Marsden, J. E.} 
\titlosb{Foundations of Mechanics}
\ekdoths{Addison-Wesley Publishing Company, Inc. (1978)}}


\Ref\bagone{
\author{Baguis, P.}
\titlosa{Induction of Hamiltonian Poisson actions}
\tenpit {Diff. Geom. Appl., \tenprm in press}}

\Ref\bagtwo{
\author{Baguis, P.}
\titlosa{Semidirect products and the Pukanszky condition}
\periodiko{J. Geom. Phys.}
\volume{25}
\selides{245--270 (1998)}}

\Ref\phd{
\author{Baguis, P.}
\titlosb{Proc\'edures de r\'eduction et d'induction en g\'eom\'etrie
symplectique et de Poisson. Applications.}
\ekdoths{Th\`ese de Doctorat, Universit\'e d'Aix-Marseille II,
D\'ecembre 1997}}

\Ref\ben{
\author{Benayed, M.}
\titlosa{Central extensions of Lie bialgebras and Poisson-Lie groups}
\periodiko{J. Geom. Phys.}
\volume{16}
\selides{301--304 (1995)}}

\Ref\drinfeld{
\author{Drinfeld, V. G.}
\titlosa{On Poisson homogeneous spaces of Poisson-Lie groups}
\periodiko{Theor. Math. Phys.}
\volume{95}
\selides{226--227 (1993)}}

\Ref\de{
\author{Duval, C., Elhadad, J.}
\titlosa{Geometric quantization and localization of
relativistic spin systems}
\periodiko{Contemp. Math.}
\volume{132}
\selides{317--330 (1992)}}

\Ref\det{
\author{Duval, C., Elhadad, J., Tuynman, G. M.}
\titlosa{Pukanszky's condition and symplectic induction}
\periodiko{J. Diff. Geom.}
\volume{36}
\selides{331--348 (1992)}}

\Ref\karolin{
\author{Karolinsky, E. A.}
\titlosa{A classification of Poisson Homogeneous spaces of complex 
reductive Poisson-Lie groups}
\titlosb{math.QA/9901073}
\selides{18 Jan 1999}}

\Ref\kazhdan{
\author{Kazhdan, D., Kostant, B., Sternberg, S.}
\titlosa{Hamiltonian Group Actions and Dynamical Systems
of Calogero Type}
\periodiko{Comm. Pure Appl. Math.}
\volume{31}
\selides{481--508 (1978)}}

\Ref\kirillov{
\author{Kirillov, A. A.}
\titlosa{Unitary representations of nilpotent Lie groups}
\periodiko{Russ. Math. Surveys}
\volume{17}
\selides{53--104 (1962)}}

\Ref\kostant{
\author{Kostant, B.}
\titlosa{Quantization and Unitary Representations}
\periodiko{Lecture Notes in Math., Springer-Verlag, Berlin}
\volume{170}
\selides{87--208 (1970)}}

\Ref\jhlu{
\author{Lu, J.-H.} 
\titlosb{Multiplicative and affine
Poisson structures on Lie Groups} 
\ekdoths{Ph.D. thesis, Univ. of California, Berkeley (1990)}}

\Ref\jhlut{
\author{Lu, J.-H.}
\titlosa{Classical dynamical r-matrices and homogeneous Poisson 
structures on $G/H$ and $K/T$}
\titlosb{math.SG/9909004 }
\selides{1 Sep 1999}}

\Ref\luwei{
\author{Lu, J.-H., Weinstein, A.}
\titlosa{Poisson-Lie groups, dressing transformations and Bruhat
decompositions}
\periodiko{J. Diff. Geom.}
\volume{31}
\selides{501--526 (1990)}}


\Ref\souriau{
\author{Souriau, J.-M.}
\titlosb{Structures des syst\`emes dynamiques}
\ekdoths{Dunod, Paris (1969)}}



\hskip10cm

\centerline{\labf Homogeneous symplectic manifolds}
\vskip0.2cm
\centerline{\labf of}
\vskip0.2cm
\centerline{\labf Poisson-Lie groups}

\vskip1cm

\centerline{\bf P. Baguis\footnote{$^{1}$}{e-mail: 
pbaguis@ulb.ac.be}$^{,}$\footnote{$^{2}$}{
Research supported by
the ``Communaut\'e fran\c caise de Belgique",
through an ``Action de Recherche Concert\'ee de la Direction de la
Recherche Scientifique".}}

\vskip0.3cm

\centerline{Universit\'e Libre de Bruxelles}
\centerline{Campus Plaine, CP 218 Bd du Triomphe}
\centerline{1050, Brussels, Belgium}

\vskip1cm

\abstract{Symplectic manifolds which are homogeneous spaces of 
Poisson-Lie groups are studied in this paper. We show that these 
spaces are, under certain assumptions, covering spaces of dressing 
orbits of the Poisson-Lie groups which act on them. The effect of the 
Poisson induction procedure on such spaces is also examined, thus leading 
to an interesting generalization of the notion of homogeneous space. 
Some examples of homogeneous spaces of Poisson-Lie groups are 
discussed in the light of the previous results.}

\vskip1cm

{\tenprm 
{\tenpit Key-words}: Poisson-Lie groups, induction of Poisson 
actions, momentum mapping,

homogeneous spaces

\vskip0.2cm

{\tenpit MSC 2000}: 53C15, 53D17, 53D20}

\vfill\eject

\baselineskip=14pt plus .2pt

\chapter{Introduction}

Homogeneous symplectic manifolds are, under reasonable conditions, 
locally isomorphic to coadjoint orbits and their relation to the 
theory of unitary irreducible
representations of Lie groups has been very early 
recognized {\kostant}. This kind of 
symplectic manifolds (together with the coadjoint orbits of Lie 
groups, which are a special case), is perhaps the most important 
non-trivial class of geometrically quantizable symplectic manifolds 
in the Kirillov-Kostant-Souriau program {\kirillov}, {\kostant}, 
{\souriau}. A fundamental ingredient of this approach is the existence 
of an equivariant momentum map for the symplectic action on the 
homogeneous symplectic manifold. Then, it turns out that such a 
symplectic manifold is a covering space of a coadjoint orbit of the 
group. This is essentially a ``Hamiltonian" classification of the 
homogeneous symplectic manifolds.

\heads

Recently, an analogous study but in a different context, has been 
initiated by {\drinfeld}, {\karolin}, {\jhlut} for Poisson-Lie groups acting on 
Poisson manifolds. In particular, a correspondence between 
Poisson homogeneous $G$-spaces, where $G$ is a Poisson-Lie group, and 
Lagrangian subalgebras of the double $D(\frak g)$ of the tangent Lie 
bialgebra of $G$, has been established in {\drinfeld}.

In the present article, we will turn our attention again to the 
Hamiltonian point of view but for Poisson-Lie groups this time. We 
first establish the exact Poisson-Lie analog of homogeneous symplectic 
manifolds: a symplectic manifold on which a Poisson-Lie group acts 
transitively through a Poisson action admitting an equivariant 
momentum map, in the sense of {\jhlu}, is a covering space of a 
dressing orbit of the Poisson-Lie group (see Proposition 2.2 below).  
Notice here that it makes no sense to replace the symplectic manifold 
by a Poisson one, because a transitive action of a Poisson-Lie group 
on a Poisson manifold has never a momentum map unless the Poisson 
manifold contains only one symplectic leaf, therefore is symplectic.

The case of non-equivariant momentum maps needs special attention in 
the Poisson-Lie case since, unlike the symplectic case, the lack 
of equivariance now is not automatically adjusted. We address this 
issue in section 3.

We also examine the effects of the Poisson induction procedure, 
introduced in {\bagone}, on a symplectic manifold which is a 
homogeneous Hamiltonian space of a Poisson-Lie group. The result 
depends on the circumstance and the procedure leads either to a 
homogeneous space or to an almost homogeneous space. The later
is introduced in Definition 2.3; actually, this notion emerges 
naturally in the induction procedure and its meaning is that the 
almost homogeneous space is generated by a discrete (eventually 
finite) subset through the action. If this subset reduces to a point, 
then we obtain a homogeneous space.

Some examples are finally discussed in section 5. More precisely, we 
describe situations in which one can have a transitive Poisson action 
on a symplectic manifold admiting a momentum map in the sense of {\jhlu} 
and we give partial solutions to the difficult problem of 
equivariance. This progressively leads to a Hamiltonian description of 
the cells of a Bruhat decomposition for coadjoint orbits of a certain type.
We finally endow a semi-direct product $G=K\timess\rho V$ with a 
Poisson-Lie structure that has the following property: the 
corresponding dressing orbits of $G$ in $G^{\ast}$ can be obtained by 
Poisson induction on coadjoint orbits of certain subgroups of $K$.

\vskip0.2cm

{\bf Conventions.} If $(P,\pi_{P})$ is a Poisson manifold, then 
$\pi_{P}^{\sharp}\colon T^{\ast}P\rightarrow TP$ is the map defined 
by 
$\alpha(\pi_{P}^{\sharp}(\beta))=\pi_{P}(\alpha,\beta),\forall\alpha,\beta\an 
T^{\ast}P$. Let now $\sigma\colon G\times 
P\rightarrow P$ (resp.  $\sigma\colon P\times 
G\rightarrow P$) be a left (resp. right) Poisson action of the Poisson-Lie group
$(G,\pi_{G})$ on $(P,\pi_{P})$,  and let us denote by $\sigma(X)$ the 
infinitesimal generator of the action and by $G^{\ast}$ the dual group 
of $G$. Then, we say that 
$\sigma$ is Hamiltonian if there exists a differentiable map $J\colon 
P\rightarrow G^{\ast}$,  called momentum mapping, 
satisfying the following equation, for each $X\an\frak g$:
$$\sigma(X)=\pi_{P}^{\sharp}(J^{\ast}X^{l})\quad(resp. \quad
\sigma(X)=-\pi_{P}^{\sharp}(J^{\ast}X^{r})).$$
In the previous equation $X^{l}$ (resp. $X^{r}$) is the left (resp. 
right) invariant 1-form on $G^{\ast}$ whose value at the identity is 
equal to $X\an\frak g\cong(\frak g^{\ast})^{\ast}$. The momentum 
mapping is said to be equivariant, if it is a morphism of Poisson 
manifolds with respect to the Poisson structure $\pi_{P}$ on $P$ and 
the canonical Poisson structure on the dual group of the Poisson Lie 
group $(G,\pi_{G})$. Left and right 
infinitesimal dressing actions $\lambda\colon\frak 
g^{\ast}\rightarrow{\eusm X}(G)$ and $\rho\colon\frak 
g^{\ast}\rightarrow{\eusm X}(G)$ of $\frak g$ on $G^{\ast}$ are defined by
$$\lambda(\xi)=\pi_{G}^{\sharp}(\xi^{l})\quad\hbox{and}\quad
\rho(\xi)=-\pi_{G}^{\sharp}(\xi^{r}),\quad\forall \xi\an\frak g^{\ast}.$$
Similarly, one defines infinitesimal left and right dressing actions 
of $\frak g$ on $G^{\ast}$. In the case where the vector fields 
$\lambda(\xi)$ (or, equivalently, $\rho(\xi)$) are complete for all 
$\xi\an\frak g^{\ast}$,  we have left and right actions of 
$(G^{\ast},\pi_{G^{\ast}})$ on $(G,\pi_{G})$ denoted also by $\lambda$ 
and $\rho$ respectively, and we say that $(G,\pi_{G})$
is a complete Poisson-Lie group.

\chapter{The equivariant Poisson-Lie case}

In the symplectic context, the following is well known. Let $(M,\omega)$ be a 
symplectic manifold and $\sigma\colon G\times M\rightarrow M$ a 
symplectic action admitting the momentum mapping $J\colon 
M\rightarrow\frak g^{\ast}$, that is the infinitesimal generator of 
the action corresponding to the element $X\an\frak g$ is equal to the 
Hamiltonian vector field corresponding to the function $J^{\ast}X\an 
C^{\infty}(M)$, where we regard the pull-back $J^{\ast}$ as a linear 
map $\frak g\rightarrow C^{\infty}(M)$.  Assume that the momentum map 
is equivariant, that is $J\comp\sigma_{g}=\Coad(g)\comp J,\forall 
g\an G$ and that the action $\sigma$ is transitive. Then:

\math{Theorem ({\kostant}).}{\shomog}{\sl Under the assumptions above, 
there exists an element $\mu_{0}\an\frak g^{\ast}$ such that $M$ be a 
covering space of the coadjoint orbit $\eusm O_{\mu_{0}}=G\cdot\mu_{0}$ and 
$J\colon M\rightarrow\eusm O_{\mu_{0}}$ be a morphism of Hamiltonian 
$G$-spaces.}

\vskip0.3cm

When the momentum map $J$ is not equivariant, one can consider the 
bilinear map $\gamma\colon\frak g\times\frak g\rightarrow C^{\infty}(M)$ 
given by $\gamma(X,Y)=\{J^{\ast}X,J^{\ast}Y\}-J^{\ast}[X,Y]$, which 
is a constant function on $M$ for each $X,Y\an\frak g$ and defines a 
2-cocycle $\frak g\times\frak g\rightarrow\R$. This completely 
determines a central extension $\tilde{\frak g}$ of the Lie algebra 
$\frak g$ and a symplectic action on $M$ of the connected and simply 
connected Lie group $\tilde{G}$ whose Lie algebra is $\tilde{\frak 
g}$. This action admits now an equivariant momentum map and the 
Theorem {\shomog} can be applied.

We consider now a symplectic manifold $P$, where the symplectic structure 
is described by a non-degenerate Poisson tensor $\pi_{P}$,  a 
Poisson-Lie group $(G,\pi_{G})$ and a left Poisson action $\sigma\colon 
G\times P\rightarrow P$. We make the assumption that the action 
$\sigma$ admits an equivariant momentum map  $J\colon P\rightarrow 
G^{\ast}$. In that case, one says that $(P,\pi_{P})$ is a {\it Hamiltonian 
Poisson-Lie $G$-space}. One has:

\math{Proposition.}{\LPhomog}{\sl If the action $\sigma$ is 
transitive, then there exists an element $u_{0}\an G^{\ast}$ such that 
the symplectic manifold $P$ be a covering space of the left dressing 
orbit $\eusm O^{l}_{u_{0}}$ of $u_{0}$ and the map $J\colon P\rightarrow
\eusm O^{l}_{u_{0}}$ be a morphism of Hamiltonian Poisson-Lie 
$G$-spaces.}

\undertext{\it Proof.} Let $p_{0}\an P$ and $u_{0}=J(p_{0})\an G^{\ast}$. 
Then, by the equivariance (and infinitesimal equivariance) of the momentum 
mapping, we find that $J\colon P\rightarrow\eusm O^{l}_{u_{0}}$ is a 
surjective submersion, where $\eusm O^{l}_{u_{0}}$ is the left 
dressing orbit of $u_{0}$. 

Consider now the equation $T_{p}J(v)=0$ for $v\an T_{p}P$. One can 
write $v=\sigma(X)_{p},X\an\frak g$ because $P$ is homogeneous space 
of $G$. Then $\lambda(X)(J(p))=0$ and consequently 
$(J^{\ast}X^{l})(\sigma(Y)_{p})=0$ for each $Y\an\frak g$. Using 
again the fact that $P$ is homogeneous and that $J\colon 
P\rightarrow\eusm O^{l}_{u_{0}}$ is a surjective submersion, we have 
the result. \qed

\vskip0.3cm

If we write $P=G/G_{p_{0}}$ and $\eusm O^{l}_{u_{0}}=G/G_{u_{0}}$, where 
$G_{p_{0}}$ and $G_{u_{0}}$ are the isotropy subgroups of $p_{0}$ and 
$u_{0}$ for the corresponding actions, then the map $J\colon 
P\rightarrow\eusm O^{l}_{u_{0}}$ can 
be writen as $J([g]_{p_{0}})=[g]_{u_{0}}$, where $[g]_{p_{0}}$ and 
$[g]_{u_{0}}$ denote the equivalence classes of $g\an G$ under the 
equivalence relations defined by the subgroups $G_{p_{0}}$ and $G_{u_{0}}$. 
Furthermore, one has $G_{p_{0}}\subset G_{u_{0}}$ and the fibre of $J$ 
is exactly $G_{u_{0}}/G_{p_{0}}$. 

There exists a generalization of the notion of homogeneous 
space which arises naturally in the induction procedure, as we will 
see later. We give the following definition:

\math{Definition.}{\almhom}{\sl Let $P$ be a differentiable manifold 
on which the Lie group $G$ acts smoothly. We will say that $P$ is an almost 
homogeneous space of $G$ if 
there exists a discrete subset $\Sigma\subset P$ such that 
$G\cdot\Sigma=P$.}

\vskip0.3cm

\noindent Otherwise stated, $P$ is almost homogeneous, if the set of $G$-orbits in 
$P$ is discrete. In particular, when $\Sigma$ is a one-point set, $P$ 
is homogeneous.

 Assume now that $(P,\pi_{P})$ is a symplectic manifold which is
an almost homogeneous space of the Poisson-Lie group $(G,\pi_{G})$
for a left Hamiltonian Poisson action of $G$ on $P$. Then, it is immediate 
from Proposition {\LPhomog} that all the open orbits of $G$ in $P$ 
are covering spaces of left dressing orbits of $G$ in $G^{\ast}$. In 
particular, if the topology on $\Sigma$, viewed as a quotient space 
$P/G$, is the discrete one, then the manifold $P$ if ``foliated'' by 
a discrete set of open submanifolds, each of them is  a covering 
space of a left dressing orbit in $G^{\ast}$.

As an example, let us  discuss the following situation coming from 
the completely symplectic setting. We take $P=\R^{2}$ equipped with 
its canonical symplectic structure and $G$ equal to the semidirect 
product between $K=\R_{+}^{\star}$ (non-zero positive real numbers)
and $V=\R$ through the 
representation $K\rightarrow GL(V)$ given by 
$r\rightarrow\displaystyle {1\over r}$ {\bagtwo}. Then,  
$(r,a)\cdot(x,y)=(rx,\displaystyle{1\over r}y+a)$ is a symplectic action 
which  
admits the equivariant momentum map $J\colon P\rightarrow \frak 
g^{\ast}\cong\R^{2}$ given by $J(x,y)=(-xy,x)$. Using the induction techniques 
of {\bagtwo} for coadjoint orbits of semidirect products, one easily 
finds that the only coadjoint orbits of $G$ are either a point or the 
cotangent bundle $T^{\ast}V\cong\R^{2}$ with its canonical symplectic 
structure. On the other hand, the space $P$ is almost homogeneous 
since it is generated by the set $\Sigma=\{(-1,0),(0,0),(1,0)\}$
through the action of $G$. Furthermore, there exist two open orbits 
in $P$, the open half-planes $P_{\pm}=\R_{\pm}^{\star}\times\R$. 
According to the previous discussion on almost homogeneous spaces, 
$P_{\pm}$ coincide with coadjoint orbits of $G$. Apparently, nothing 
can be said about the $G$-orbit of $(0,0)$ which coincides with 
the $y$-axis and therfore is 1-dimensional.

\chapter{The non-equivariant Poisson-Lie case}

We are now placed in the case where the Poisson action $\sigma\colon 
G\times P\rightarrow P$ admits a non-equivariant momentum map 
$J\colon P\rightarrow G^{\ast}$. Without loss of generality, we can 
assume that there exists a point $x_{0}\an P$ such that 
$J(x_{0})=e^{\ast}$, the identity of $G^{\ast}$;
let $\gamma=J_{\ast}\pi_{P}(x_{0})\an\Lambda^{2}\frak g^{\ast}$.
For each $X,Y\an\frak g$, we 
consider the function $\mu(X,Y)\an C^{\infty}(P)$ defined by
$$\mu(X,Y)=\pi_{P}(J^{\ast}X^{l},J^{\ast}Y^{l})-
J^{\ast}(\pi_{G^{\ast}}(X^{l},Y^{l})).\eqn\muxy$$
The function $\mu(X,Y)$ controls the equivariance of $J$ because 
$\mu(X,Y)=0$ forall $X,Y\an\frak g$ if and only if $J$ is equivariant. 

Let now $\pi_{J}=\pi_{G^{\ast}}+\gamma^{r}$. Then, it can be proved 
{\jhlu} that $J\colon(P,\pi_{P})\rightarrow(G^{\ast},\pi_{J})$ is a Poisson 
map. But this means that
$$\mu(X,Y)(p)=(\Ad(J(p)^{-1})\gamma)(X,Y),\forall X,Y\an\frak g,p\an 
P.$$
We see that if $\gamma$ is Ad-invariant, then $\mu(X,Y)$ is a 
constant function on $P$. But there is more than this:

\math{Proposition.}{\gammacocycle}{\sl If $\gamma$ is Ad-invariant, 
then $\gamma$ is a real-valued 2-cocycle on the Lie algebra $\frak g$.}

\vskip0.3cm

\noindent We refer the reader to {\phd} for the proof of this proposition.
Assuming now that $\gamma$ is Ad-invariant, let us consider the 
central extension of the Lie bialgebra $\frak g$ defined by the 
cocycle $\gamma$ and the zero derivation on $\frak g^{\ast}$. For the 
reader's convenience, we recall from {\ben} that a central extension 
of a Lie bialgebra $\frak g$ is an exact sequence 
$0\longrightarrow\R\buildrel i\over\longrightarrow \hat{\frak g}
\buildrel j\over\longrightarrow\frak g\longrightarrow 0$, such that 
$i$ and $j$ be morphisms of Lie bialgebras and $i(\R)$ be contained 
in the center of $\hat{\frak g}$. If $Z^{2}(\frak g,\R)$ and $\eusm 
Der(\frak g^{\ast})$ stand for the space of real-valued 2-cocycles 
over $\frak g$ and the space of derivations on the Lie algebra $\frak 
g^{\ast}$, then such an extension is completely determined 
by a pair of elements $(\gamma,f)\an Z^{2}(\frak g,\R)\times\eusm Der(\frak 
g^{\ast})$ which are Drinfeld-compatible:
$$f^{\ast}([X,Y])-[f^{\ast}(X),Y]-[X,f^{\ast}(Y)]=\coad(\gamma^{\sharp}(X))(Y)-
\coad(\gamma^{\sharp}(Y))(X).\eqn\drinfcomp$$

Returning to our case, the pair $(\gamma=J_{\ast}\pi_{P}(x_{0}),f=0)$ 
defines actually a central extension of $\frak g$ because the 
Ad-invariance of $\gamma$ gives the Drinfeld-compatibility with $f=0$.
If $\hat{\frak g}$ is the extended Lie bialgebra, we have an  
infinitesimal left action $\hat{\sigma}\colon\hat{\frak g}\rightarrow\eusm 
X(P)$ given by $\hat{\sigma}(X,a)=\sigma(X)$. This action is Poisson 
giving thus rise to a Poisson action, still denoted by 
$\hat{\sigma}$, of the connected and simply connected Poisson-Lie 
group $\hat{G}$ having $\hat{\frak g}$ as tangent Lie bialgebra, on the Poisson 
manifold $(P,\pi_{P})$. Under the isomorphism 
$\hat{G}\cong\tilde{G}\times\R$, where $\tilde{G}$ is the universal 
covering space of $G$, one has:

\math{Theorem.}{\nonequiv}{\sl The differentiable map $\hat{J}\colon 
P\rightarrow\hat{G}^{\ast}$ given by
$$\hat{J}(p)=(J(p),1),\forall p\an P\eqn\jhat$$
is an equivariant momentum map for the Poisson action $\hat{\sigma}$.}

\undertext{\it Proof.} We first note that the Poisson structure on 
$\hat{G}^{\ast}\cong G^{\ast}\times\R$ is given by
$$\pi_{\hat{G}^{\ast}}(v,a)=\pi_{G^{\ast}}(v)+a\gamma^{r}(v),
\forall(v,a)\an\hat{G}^{\ast}.$$
The rest of the proof consists of a calculation of the function 
$\mu(X,Y)$ for the momentum map $\hat{J}$ and it is ommited here.\qed

\chapter{Poisson induction of homogeneous symplectic manifolds}

Symplectic induction is a procedure by means of which one can induce, 
from precise initial data, symplectic structures and Hamiltonian 
actions to bigger manifolds and has several interesting applications
{\kazhdan}, {\de}, {\det}, {\bagtwo}. 
We know {\phd} that the induction procedure has a natural Poisson 
analog and, in particular, that the Poisson induced of a symplectic 
manifold is also a symplectic manifold carrying a Hamiltonian Poisson 
action. Here we will study the effect of the Poisson induction 
procedure on a homogeneous symplectic manifold, a problem whose even 
the completely symplectic aspects are not yet known.

Let us recall from {\phd} the elements of the Poisson induction which are 
necessary for the understanding of what follows. We consider a left Poisson 
action $\sigma\colon 
(H,\pi_{H})\times(P,\pi_{P})\rightarrow(P,\pi_{P})$ admitting the 
equivariant momentum mapping $J\colon P\rightarrow H^{\ast}$ and 
assume that $(H,\pi_{H})$ is a Poisson-Lie subgroup of the Poisson-Lie 
group $(G,\pi_{G})$. In order to simplify technically the discussion, 
we make the assumption that $(G,\pi_{G})$ is simply connected and 
complete. If $(D(G),\pi_{+})$ is the double group of $G$ equipped with 
its natural symplectic structure {\jhlu}, let 
$(\tilde{P},\pi_{\tilde{P}})=(P,\pi_{P})\times(D(G),\pi_{+})$. Then 
the map $\tilde{\sigma}\colon H\times\tilde{P}\rightarrow\tilde{P}$ 
given by
$$\tilde{\sigma}(h,(p,d))=(\sigma(\lambda_{i^{\ast}u^{-1}}(h),p),
d[\lambda_{i^{\ast}u^{-1}}(h)]^{-1}),\forall h\an H,(p,d)\an\tilde{P}
\eqn\sigmatilde$$
where $i^{\ast}\colon G^{\ast}\rightarrow H^{\ast}$ is the projection 
of the dual groups induced by the inclusion $i\colon H\hookrightarrow 
G$, is a left Poisson action of $H$ on $\tilde{P}$ admitting the 
equivariant momentum map $\tilde{J}\colon\tilde{P}\rightarrow 
H^{\ast}$ with 
$$\tilde{J}(p,d)=J(p)(i^{\ast}u^{-1}).\eqn\jtilde$$
If $e^{\ast}\an H^{\ast}$ is the unit element, we obtain by
Marsden-Weinstein reduction the induced manifold as 
$$P_{ind}=\tilde{J}^{-1}(e^{\ast})/H.\eqn\pind$$
The group $G$ acts on $P_{ind}$ as follows: if $[(p,gu)]$ is the 
equivalence class of $(p,gu)\an\tilde{J}^{-1}(e^{\ast})$ in the 
quotient $\pind$, then 
$$k\cdot[(p,gu)]=[(p,\lambda_{\rho_{g^{-1}}(u)}(k)gu)]\eqn\gaction$$
for all $k\an 
G$. This action is Poisson and admits an equivariant momenutm mapping.

The above construction works even when $(P,\pi_{P})$ is a general 
Poisson manifold, but here we are interested in the case where $P$ is 
symplectic. Then, $P_{ind}$ is also symplectic and the 
following question arises naturally: what can we tell about 
$P_{ind}$ when $P$ is a homogeneous space of $H$?

Let us fix an element $[(p_{0},g_{0}u_{0})]\an P_{ind}$ and consider 
an arbitrary element $[(p,gu)]\an P_{ind}$. Then, assuming that $H$ 
acts transitively on $P$, we conclude that there exist elements $h\an 
H$ and $u^{1}\an H^{\circ}$ such that
$$\left\{\matrix{p&=&(\sigma_{J})_{h}(p_{0})\cr
u&=&\rho_{h^{-1}}(u_{0}u^{1}).\cr}\right.\eqn\eqone$$
Here $H^{\circ}=\ker(i^{\ast})$ and $\sigma_{J}$ is the action of $H$ 
on $P$ defined as 
$(\sigma_{J})_{h}(p)=\sigma_{\lambda_{J(p)^{-1}}(h)}(p)$.
We want now to solve the equation 
$k\cdot[(p_{0},g_{0}u_{0})]=[(p,gu)]$ with respect to $k\an G$ for 
given $[(p_{0},g_{0}u_{0})]$ and arbitrary $[(p,gu)]$. In view of the 
transitivity of the action $\sigma$ and of $\gaction$, this equation is equivalent to
$$\left\{\matrix{p&=&(\sigma_{J})_{h}(p_{0})\cr
g&=&\lambda_{\rho_{g_{0}^{-1}}(u_{0})}(k)g_{0}h^{-1}\cr
u&=&\rho_{h^{-1}}(u_{0}).\cr}\right.\eqn\eqtwo$$
The last equations make clear that, generally, one can nothing say about 
the action of $G$ on $P_{ind}$, despite the transitivity of $\sigma$. 
Let us make the following assumption: $u_{0}H^{\circ}\subset H\cdot 
u_{0}$, where the dot on the right hand side means left dressing 
transformations. If $J(p_{0})=v_{0}$, one then observes:

\item{$\bullet$} if $h_{1}\an H$ is an element for which 
$\lambda_{h_{1}}(u_{0})=u_{0}u^{1}$, then $h_{1}\an H_{v_{0}}$, 
where $H_{v_{0}}$ represents the isotropy subgroup of $v_{0}\an 
H^{\ast}$ with respect to the left dressing tranformations;
\item{$\bullet$} $u_{0}u^{1}=\rho_{\tilde{h}_{1}^{-1}}(u_{0})$, 
where $\tilde{h}_{1}^{-1}=\lambda_{J(p_{0})}(h_{1})$.

\noindent Consequently, the equations $\eqone$ become
$$\left\{\matrix{p&=&(\sigma_{J})_{h\tilde{h}_{1}}(p_{0}^{\prime})\cr
u&=&\rho_{(h\tilde{h}_{1})^{-1}}(u_{0}),\cr}\right.$$
where $p_{0}^{\prime}=\sigma_{h_{1}^{-1}}(p_{0})\an J^{-1}(v_{0})$.
Combining with equations $\eqtwo$ we obtain the following theorem:

\math{Theorem.}{\homindution}{\sl Assume that the symplectic 
manifold $(P,\pi_{P})$ is a homogeneous Hamiltonian Poisson-Lie space of 
the Poisson-Lie group $(H,\pi_{H})$, viewed as a Poisson-Lie subgroup 
of $(G,\pi_{G})$. Then, the Poisson induced symplectic manifold 
$(P_{ind},\pi_{ind})$ has the following structure:
\item{$\bullet$} if $u_{0}H^{\circ}\subset H\cdot u_{0}$, then 
$P_{ind}$ is an almost homogeneous Hamiltonian Poisson-Lie space of $G$;
\item{$\bullet$} if $u_{0}H^{\circ}\subset H_{p_{0}}\cdot u_{0}$, 
where $H_{p_{0}}$ is the isotropy subgroup of $p_{0}$ for the 
action $\sigma$ of $H$ on $P$, then $P_{ind}$ is a homogeneous  
Hamiltonian Poisson-Lie space of $G$ and, in view of Proposition
{\LPhomog}, a covering space of a left dressing orbit of $G$ in $G^{\ast}$.}

\chapter{Applications}

{\bf (1)} Consider two complete Poisson-Lie groups $(H,\pi_{H})$ and 
$(G,\pi_{G})$ and an injective morphism of Poisson-Lie groups 
$f\colon H\rightarrow G$ which is an immersion at the identity. Then, 
the linear map $\sigma\colon\frak g^{\ast}
\rightarrow\eusm X(H)$ given by
$$\sigma(\xi)=\pi_{H}^{\sharp}(f^{\ast}\xi^{l})$$
is an infinitesimal Poisson action admitting $f$ as equivariant 
momentum map. In the complete case we are studying, the corresponding
Poisson action $\sigma$ of $G^{\ast}$ on $H$ is given by the equation
$\sigma_{u}(h)=\lambda_{f^{\ast}u}(h)$, where $f^{\ast}\colon 
G^{\ast}\rightarrow H^{\ast}$ is the morphism of the dual groups 
induced by $f$. In view of Proposition {\LPhomog}, one obtains:

\math{Corollary.}{\HGorbits}{\sl 
The orbit of $h\an H$ under the left dressing 
action of $H^{\ast}$ is
a homogeneous Hamiltonian Poisson-Lie 
$G^{\ast}$-space and, at the same time,
a covering space of the orbit of $f(h)\an G$ 
under the left dressing action of $G^{\ast}$. In particular, these 
orbits have the same dimension.}

\vskip0.3cm

{\bf (2)} If $(P,\pi_{P})$ is a non-trivial 
homogeneous Hamiltonian Poisson-Lie space 
of $(G,\pi_{G})$, then there is no point of $P$ whose image under the 
momentum map $J$ could be the identity of the group $G^{\ast}$.
Indeed, if such a point existed, then $P$ should be locally isomorphic 
to a point according to Proposition {\LPhomog}, which is a 
contradiction.

Consider now a simply connected symplectic manifold  $P$ 
and a transitive
left Poisson action $\sigma\colon G\times P\rightarrow P$ of 
the Poisson-Lie group $(G,\pi_{G})$. Then, for each point $x_{0}\an 
P$ there exists a momentum 
mapping $J\colon P\rightarrow G^{\ast}$ for this action {\jhlu}, such that 
$J(x_{0})=e^{\ast}$, the identity of $G^{\ast}$. Then, according to 
the previous argument, this momentum map cannot be equivariant.
Assume instead that there exists 
an element $u_{0}\an G^{\ast}$, such that the bivectors $\pi_{P}$ and 
$\pi_{G^{\ast}}$ be $J_{0}$-related at the point $x_{0}$, where
$J_{0}=L_{u_{0}}\comp J$. Then, $J_{0}$ is a Poisson morphism 
between $(P,\pi_{P})$ and $(G,\pi_{G^{\ast}})$ and hence an 
equivariant momentum map for the action $\sigma$. In this case, $P$ is 
the universal covering space of a left dressing orbit of $G$ in $G^{\ast}$.

Staying always in the case of a transitive left Poisson action on a 
simply connected symplectic manifold, assume that there exists a point 
$x_{0}\an P$ such that the annihilator of the isotropy subalgebra
$\frak g_{x_{0}}$ (with respect to the action $\sigma$)
be contained in the center of the Lie algebra $\frak 
g^{\ast}$. Then, the  element $\gamma\an\Lambda^{2}\frak g$
of Proposition {\gammacocycle} is invariant under the adjoint
action of $G^{\ast}$ but not zero, because of the transitivity of the 
action. Then, by Theorem {\nonequiv} and Proposition {\LPhomog},
$P$ is the universal covering 
space of a left dressing orbit in the dual group of an appropriate 
central extension of $G$.

\vskip0.3cm

{\bf (3)} We consider now the case where $(H,\pi_{H})$ is a Poisson-Lie 
group of $(G,\pi_{G})$, with $\pi_{G}(g)=R_{g}r-L_{g}r$ and 
$\pi_{H}=0$, where $r\an\Lambda^{2}\frak g$ has the property 
$\Ad(g)[r,r]=[r,r]$. Otherwise stated, $\pi_{G}$ is exact and 
$\Ad(h)r=r,\forall h\an H$. Let $\eusm O=G/H$ and 
$\bar{\rho}\colon\eusm O\times G^{\ast}\rightarrow\eusm O$ the right 
action of $G^{\ast}$ on the Poisson manifold $\eusm O$ {\luwei}, 
obtained projecting the right dressing transformations of $G^{\ast}$ 
in $G$.   
Then, for each point $x_{0}=[g_{0}]\an\eusm O$, the orbit $x_{0}\cdot 
G^{\ast}$ coincides with the symplectic leaf through $x_{0}$. Pick an 
element $g_{0}\an G$ such that the orbit of $x_{0}$ be simply 
connected. Then, the action $\bar{\rho}$ restricts to a right Poisson 
action of $G^{\ast}$ on $x_{0}\cdot G^{\ast}$ for which
there exists a momentum mapping $J\colon x_{0}\cdot 
G^{\ast}\rightarrow G$ with $J(x_{0})=e$. By construction, the tangent 
of $J$ at $x_{0}$ is given by $J_{\ast 
x_{0}}(\bar{\rho}(\eta)_{x_{0}})(\xi)=-\omega_{x_{0}}(\bar{\rho}
(\eta)_{x_{0}},\bar{\rho}(\xi)_{x_{0}}),\forall\eta,\xi\an\frak 
g^{\ast}$, where $\omega$ is the symplectic structure of the 
symplectic leaf $x_{0}\cdot G^{\ast}$. 
Assume now that the Lie algebra $\frak h$ of $H$ has center of 
dimension 0 or 1. Then, one finds
$$J_{\ast x_{0}}(\bar{\rho}(\eta)_{x_{0}})=R_{g_{0}^{-1}}
\rho(\eta)_{g_{0}}$$
which has as direct consequence that $J_{\ast 
x_{0}}\pi(x_{0})=R_{g_{0}^{-1}}\pi_{G}(g_{0})$, where $\pi$ is the 
Poisson structure of $\eusm O$. Then, the results of 
example (2) above, conveniently adapted for a right Poisson action, 
confirm that $J_{0}=R_{g_{0}}\comp J\colon x_{0}\cdot 
G^{\ast}\rightarrow G$ is an equivariant momentum map for the 
transitive Poisson action $\bar{\rho}$. Consequently:

\math{Proposition.}{\sconleaf}{\sl Let $(G,\pi_{G})$ be a Poisson-Lie 
group, where $\pi_{G}$ is exact with linearization at the identity 
equal to $r\an\Lambda^{2}\frak g$.
Assume further that $r$ is invariant under the 
adjoint action of a closed Lie subgroup $H$ and that the center of $\frak 
h=Lie(H)$
is equal to $\R X_{0}$, where $X_{0}\an\frak h$ is eventually zero.
Then, all the 
simply connected symplectic leaves of the quotient Poisson space $G/H$ 
are universal covering spaces of right dressing orbits of $G^{\ast}$ 
in $G$.}

\vskip0.3cm

In particular, when $G=SU(2)$ viewed as a Poisson-Lie group as in 
{\luwei} and $H=\S^{1}$ with the zero Poisson structure, the Poisson 
manifold $\eusm O$ coincides with the 2-sphere $\S^{2}$. In this 
case, we have two symplectic leaves, a point and its complement 
(isomorphic to the plane), and the assumptions of Proposition 
{\sconleaf} are fulfilled. Consequently, the open leaf is the 
universal covering space of a right dressing orbit of the three 
dimensional "book" group (the dual group of $SU(2)$) in $SU(2)$.
More generally, if $G$ is a compact semisimple Lie group and $\eusm O$ a 
coadjoint orbit of $G$, then according to {\luwei}, $G$ can be equipped 
with a Poisson-Lie structure which descends to $\eusm O$ as a Poisson 
structure whose symplectic leaves are all cells of a Bruhat 
decomposition of $\eusm O$ and diffeomorphic to dressing orbits of 
$G^{\ast}$ in $G$. The Propositions {\LPhomog} and {\sconleaf} provide 
a Hamiltonian interpretation of this situation and indicate a partial 
generalization for arbitrary Poisson-Lie group $G$.

\vskip0.3cm

{\bf (4)} We discuss now an example where the conditions of Theorem 
{\homindution} are always satisfied.
Let us consider a semi-direct product $G=K\timess\rho \frak v$, 
formed by a Lie group $K$ and a vector space $\frak v$ through the 
representation $\rho\colon K\rightarrow GL(\frak v)$. 
The group law on $G$ is given by
$$(\kappa,u)\cdot(\lambda,v)=(\kappa\lambda,\kappa\cdot v+u),
\eqn\semlaw$$
for all $(\kappa,u),(\lambda,v)\an G$, where $\kappa\cdot v$ means the 
action of the element $\kappa\an K$ on $v\an \frak v$ through the 
representation $\rho$ {\bagtwo}. If $\phi\colon 
\frak v\rightarrow\Lambda^{2}T\frak v$ 
is a Poisson-Lie structure on $\frak v$ (for the 
abelian group law of a vector space) invariant under the 
representation $\rho$, then  the Poisson tensor 
$\pi=0\oplus\phi$ is a Poisson-Lie structure on $G$ for the group 
operation given by $\semlaw$. Now, $\frak v^{\ast}$ inherits a Lie algebra 
structure and we will denote by $V^{\ast}$ the corresponding connected and 
simply connected Lie group with Lie algebra $\frak v^{\ast}$. The Lie 
group $V^{\ast}$ can obviously be seen as the dual group of the 
Poisson-Lie group $V=(\frak v,\phi)$. If $\frak k$ is the Lie algebra 
of $K$, then the dual group of $G$ is
$$G^{\ast}=\frak k^{\ast}\times V^{\ast}$$
equipped with the direct product group operation. For $X\an\frak k$ 
and $\kappa\an K$, one has two mappings:
$$\frak v^{\ast}\rightarrow\frak 
v^{\ast},\quad q\mapsto X\cdot q,\forall q
\an\frak v^{\ast}\eqn\mux$$
and
$$\frak v^{\ast}\rightarrow\frak 
v^{\ast},\quad q\mapsto\kappa\cdot q,\forall q\an\frak 
v^{\ast}.\eqn\muk$$
In the previous equations, $X\cdot q$ and $\kappa\cdot q$ mean the 
actions of $\frak k$ and $K$ respectively through the contragredient 
representation on $\frak v^{\ast}$. The map $\muk$ can be 
integrated to an action (by group homomorphisms) 
of the group $K$ on the dual group $V^{\ast}$. The induced fundamental 
vector fields coincide with the multiplicative vector fields obtained 
from $\mux$ (observe that $q\mapsto X\cdot q$ is a 1-cocycle for the adjoint 
representation of $\frak v^{\ast}$ on itself). If $\varrho\colon 
K\times V^{\ast}\rightarrow V^{\ast}$ denotes this action, then one 
obtains a linear map  $\tau_{p}\colon \frak k\rightarrow\frak 
v^{\ast}$, for each $p\an V^{\ast}$, as follows:
$$\tau_{p}(X)=R_{p^{-1}}\varrho(X)_{p}.\eqn\nltau$$
We introduce the notation $\tau_{p}^{\ast}(v)=p\odot v$ for the 
transposed map $\tau_{p}^{\ast}\colon \frak v\rightarrow\frak 
k^{\ast}$. With the above data one can calculate the left dressing 
tranformation of $G$ on $G^{\ast}$ and the result is
$$\lambda_{(\kappa,a)}(\xi,p)=(\Coad(\kappa)\xi+\varrho_{\kappa}
(p)\odot a,\varrho_{\kappa}(p)),\quad\forall (\kappa,a)
\an G,(\xi,p)\an G^{\ast}
.\eqn\semleftdress$$
For a given element $u_{0}=(\xi_{0},p_{0})\an G^{\ast}$, it is 
easy to make the 
following observations using equation $\semleftdress$:

\item{$\bullet$} If $K_{p_{0}}$ is the isotropy subgroup of $p_{0}$ 
with respect to the action $\varrho$ of $K$ on $V^{\ast}$, $\frak 
k_{p_{0}}$ its Lie alebra and $i_{p_{0}}\colon\frak 
k_{p_{0}}\hookrightarrow \frak k$ the natural inclusion, then
$$\lambda_{(\kappa, a)}(i_{p_{0}}^{\ast}\xi,p_{0})=(\Coad(\kappa)
i_{p_{0}}^{\ast}\xi,p_{0}),\quad\forall (\kappa, a)\an 
K_{p_{0}}\timess\rho\frak v.$$
\item{$\bullet$} If $H=K_{p_{0}}\timess\rho\frak v$, then $H^{\circ}=
(\frak k_{p_{0}})^{\circ}\times\{e\}$ and
$$u_{0}\cdot H^{\circ}=(\xi_{0}+(\frak k_{p_{0}})^{\circ},p_{0}).$$
\item{$\bullet$} If $P$ is the orbit of $v_{0}=i^{\ast}u_{0}=
(i_{p_{0}}^{\ast}\xi,p_{0})\an H^{\ast}$ 
under the left dressing transformations of $H$, and $H_{v_{0}}$ the 
isotropy subgroup of $v_{0}$ with respect to the left dressing 
transformations of $H$ on $H^{\ast}$ , then
$$u_{0}\cdot H^{\circ}=H_{v_{0}}\cdot u_{0}.$$

\noindent Now, in the case where $P$ is a dressing orbit in 
Theorem {\homindution},
we observe that $P_{ind}$ is a dressing orbit too (see also {\bagone}).
Combining this 
remark with the previous observations, we conclude that the dressing 
orbits of $G=K\timess\rho\frak v$ (equipped with the Poisson-Lie 
structure $\pi=0\oplus\phi$) are obtained by Poisson induction on 
{\it coadjoint} orbits of certain subgroups of $K$.

Of particular interest is the case where $\frak v=\frak k^{\ast}$ and 
$\rho=\Coad$. Then $G=T^{\ast}K$ and $\pi$ 
coincides with the canonical Poisson structure on $T^{\ast}K$.
Consequently, the  dressing orbits of $T^{\ast}K$
can be obtained by Poisson induction on coadjoint orbits of subgroups 
of $K$. 

\vskip0.2cm

{\bf Acknowledgments.} I would like to thank J.-H. Lu for 
careful reading of the manuscript and critical and fruitful remarks.
It is also a pleasure to thank M. Cahen and C. Duval for their 
kind interest in this work and useful suggestions.

\vfill\eject
\vskip1cm
   \ifreferenceopen \Closeout\referencewrite \referenceopenfalse \fi
   \line{\bf\hskip0pt\hfil References\hfil}\vskip\headskip
   \vskip0.3cm
   \input referenc.txa

\end